\documentstyle[12pt]{article}
\input amssymb.sty
\newcommand{\G}{\Gamma}
\newcommand{\be}{\beta}
\newcommand{\te}{\theta}
\newcommand{\la}{\lambda}
\newcommand{\La}{\Lambda}

\newcommand{\al}{\alpha}
\newcommand{\de}{\delta}

\newtheorem{predl}{Proposition}[section]

\newtheorem{cor}{Corollary}[section]
\newtheorem{lem}{Lemma}[section]

\newcommand{\beq}[1]{\begin{equation}\label{#1}}
\newcommand{\eq}{\end{equation}}

\begin{document}
\begin{flushright}
 MIIT-TH/01-06
\end{flushright}
\vspace{10mm}
\begin{center}
{\Large \bf The asymptotic behavior of $q$-exponentials and
$q^2$-Bessel functions}\\
\vspace{5mm}
V.-B.K.Rogov \\
{\sf MIIT, 127994, Moscow, Russia} \\
{\em vrogov@cemi.rssi.ru}\\
\vspace{5mm}
\end{center}

\vspace{10mm}
\begin{abstract}
The connections between $q^2$-Bessel functions of three types and q-exponential
of three types are established. The q-exponentials and the $q^2$-Bessel functions
are represented as the Laurent series. The asymptotic behaviour of  the
q-exponentials and the $q^2$-Bessel functions is investigated.
\end{abstract}

\section{Introduction}
\setcounter{equation}{0}

The main goal of this paper is an investigation of the asymptotic behavior
of the $q^2$-Bessel functions for the large value of the argument. To solve
this problem we establish  connections between the $q$-exponentials and the
$q^2$-Bessel functions.

The $q$-exponentials of type 1 and 2 are well-known. Although the third
$q$-exponential is not very familiar nevertheless its properties can be
received from the properties of first and second $q$-exponentials.
The situation with the $q^2$-Bessel function of type 3 (the Hahn-Exton
function) is similar.

Let $q$ be the real number from $(0,1)$.
The $q$-exponentials and the $q^2$-Bessel functions are the solutions of the
difference equations, and so to consider them on the $q$-lattice (or
$q^2$-lattice) is naturally. On the other hand these functions are
determined as the convergent series in some region, and so they are analytic
ones in the corresponding regions. So we will consider the continuous
functions of $z>0$ and $z=q^{n+\nu},$ $ ~n\in\mathbb{Z},$ $ ~0\le\nu<1$. We
will investigate the asymptotic behavior of functions $f(z)$ if $z\to\infty$
or $z=q^{n+\nu}$ and $n\to-\infty$.

The main result is contained in the Propositions \ref{p4.4}, \ref{p6.1} and
in the formulas (\ref{7.4}) - (\ref{7.7}).

We will use the standard notations
$$
(a,q)_n=(1-a)(1-aq)...(1-aq^{n-1}),
~~~(a,q)_\infty=\lim_{n\to\infty}(a,q)_n,
$$
$$
(a_1,...,a_k;q)_n=(a_1,q)_n...(a_k,q)_n,
$$
and the basic hypergeometric series
\beq{1.1}
\phantom._r\Phi_s(a_1,...,a_r;b_1,...,b_s;q;z)=
\sum_{n=0}^\infty\frac{(a_1,...,a_r;q)_n}{(q,q)_n(b_1,...,b_s;q)_n}
((-1)^nq^{\frac{n(n-1)}{2}})^{s-r+1}z^n
\eq

Let $u=(1-q^2)z$. There are three types of the $q^2$-exponentials:
\beq{1.2}
1. ~e_q^{(1)}(u)=e_q((1-q^2)z)=\frac1{((1-q^2)z,q)_\infty}=
\sum_{n=0}^\infty\frac{(1-q^2)^nz^n}{(q,q)_n}, ~~~~|z|<\frac1{1-q^2},
\eq
\beq{1.3}
2. ~~e_q^{(2)}(u)=E_q((1-q^2)z)=(-(1-q^2)z,q)_\infty=
\sum_{n=0}^\infty\frac{q^{\frac{n(n-1)}2}(1-q^2)^nz^n}{(q,q)_n},
\eq
\beq{1.4}
3. ~~e_q^{(3)}(u)=\phantom._1\Phi_1(0;-\sqrt q;\sqrt q,-(1-q^2)z)=
\sum_{n=0}^\infty \frac{q^{\frac{n(n-1)}4}(1-q^2)^nz^n}{(q,q)_n}.
\eq

Obviously $e_q^{(1)}(u)$ is a meromorphic function and has simple poles at
the points $u=q^{-n}, ~~n=0, 1, ...$. Consequently,
$e_q^{(2)}(u)$ and $e_q^{(3)}(u)$ are the holomorphic functions.

Remark that for any $j=1, ~2, ~3$
$$
\lim_{q\to 1-0}e_q^{(j)}((1-q^2)z)=e^{2z}.
$$

Let
$$
D_zf(z)=\frac{f(z)-f(qz)}{(1-q^2)z}.
$$
The $q$-exponentials satisfy the next difference equations:
$$
D_ze_q^{(1)}(u)=e_q^{(1)}(u),
$$
$$
D_ze_q^{(2)}(u)=e_q^{(2)}(qu),
$$
$$
D_ze_q^{(3)}(u)=e_q^{(3)}(q^{\frac12}u).
$$

The $q-\G$-function is determined as
$$
\G_q(\al)=\frac{(q,q)_\infty}{(q^\al,q)_\infty}(1-q)^{1-\al}.
$$

\section{The $q^2$-Bessel functions }
\setcounter{equation}{0}

Remind that the $q^2$-Bessel functions of type $j, ~j=1, 2, 3,$ are determined
by the following series \cite{GR}
$$
J_\nu^{(j)}(2(1-q^2)z;q^2)=\frac1{\G_{q^2}(\nu+1)}\sum_{n=0}^\infty(-1)^n
\frac{q^{(2-\de)n(n+\nu)}(1-q^2)^{2n}z^{\nu+2n}}{(q^2,q^2)_n(q^{2\nu+2},q^2)_n},
$$
where the parameters $j$ and $\de$ are connected by formula
$$
j=-\frac32\de^2+\frac52\de+2.
$$
For $j=1, ~|z|<\frac1{1-q^2}$.
The modified $q^2$-Bessel functions are
\beq{2.1}
I_\nu^{(j)}(2(1-q^2)z;q^2)=\frac1{\G_{q^2}(\nu+1)}\sum_{n=0}^\infty
\frac{q^{(2-\de)n(n+\nu)}(1-q^2)^{2n}z^{\nu+2n}}{(q^2,q^2)_n(q^{2\nu+2},q^2)_n}.
\eq
Obviously
\beq{2.2}
I_\nu^{(j)}(ze^{-i\frac\pi2};q^2)=e^{-i\frac\nu2\pi}J_\nu^{(j)}(z;q^2),
~~~I_\nu^{(j)}(ze^{i\frac\pi2};q^2)=e^{i\frac\nu2\pi}J_\nu^{(j)}(z;q^2),
\eq

If $\nu$ is not integer the $q^2$-Neumann functions and the $q^2$-Macdonald
functions respectively
$$
Y_\nu^{(j)}(2(1-q^2)z;q^2)=\frac{q^{-\nu^2+\nu}}{\pi}\G_{q^2}(\nu)\G_{q^2}(1-\nu)
[\cos\nu\pi J_\nu^{(j)}(2(1-q^2)z;q^2)-J_{-\nu}^{(j)}(2(1-q^2)z;q^2)],
$$
$$
K_\nu^{(j)}(2(1-q^2)z;q^2)=\frac{q^{-\nu^2+\nu}}2\G_{q^2}(\nu)\G_{q^2}(1-\nu)
[I_{-\nu}^{(j)}(2(1-q^2)z;q^2)-I_\nu^{(j)}(2(1-q^2)z;q^2)].
$$

$q^2$-Wronskian of two solutions of a difference equation of second order is
determined as
$$
W(f_1,f_2)=f_1(z)f_2(qz)-f_1(qz)f_2(z).
$$

\begin{predl}\label{p2.1}
The functions $J_\nu^{(j)}(2(1-q^2)z;q^2)$ and $Y_\nu^{(j)}(2(1-q^2)z;q^2)$
satisfy the difference equation
\beq{2.3}
f(q^{-1}z)-(q^{-\nu}+q^\nu)f(z)+f(qz)=-q^{-\de}(1-q^2)^2z^2f(q^{1-\de}z),
\eq
and form the fundamental system of the solutions of this equation with
$q^2$-Wronskian
$$
W(J_\nu^{(j)},Y_\nu^{(j)})=\left\{\begin{array}{lcl}
\frac{q^{-\nu^2}(1-q^2)}\pi e_{q^2}(-(1-q^2)^2z^2)&for&\de=2\\
\frac{q^{-\nu^2}(1-q^2)}\pi&for&\de=1\\
\frac{q^{-\nu^2}(1-q^2)}\pi E_{q^2}((1-q^2)^2z^2)&for&\de=0.\\
\end{array}\right.
$$
\end{predl}
{\it Proof.} The first statement is checked directly.

Now consider the Wronskian
$$
W(J_\nu^{(j)}((1-q^2)z;q^2),Y_\nu^{(j)}((1-q^2)z;q^2))=W(J_\nu^{(j)},Y_\nu^{(j)})(z)=
$$
$$
J_\nu^{(j)}((1-q^2)z;q^2)Y_\nu^{(j)}((1-q^2)qz;q^2)-
J_\nu^{(j)}((1-q^2)qz;q^2)Y_\nu^{(j)}((1-q^2)z;q^2).
$$
Let $\de=2 (j=1)$. Because the functions $J_\nu^{(1)}$ and $Y_\nu^{(1)}$ satisfy
to the equation (\ref{2.3}) we have
$$
J_\nu^{(1)}((1-q^2)q^2z;q^2)=-(1+(1-q^2)^2z^2)J_\nu^{(1)}((1-q^2)z;q^2)+
(q^{-\nu}+q^\nu)J_\nu^{(1)}((1-q^2)qz;q^2)
$$
and
$$
Y_\nu^{(1)}((1-q^2)q^2z;q^2)=-(1+(1-q^2)^2z^2)Y_\nu^{(1)}((1-q^2)z;q^2)+
(q^{-\nu}+q^\nu)Y_\nu^{(1)}((1-q^2)qz;q^2).
$$
Hence
$$
W(J_\nu^{(1)},Y_\nu^{(1)})(qz)=(1+(1-q^2)^2z^2)W(J_\nu^{(1)},Y_\nu^{(1)})(z),
$$
i.e.
$$
W(J_\nu^{(1)},Y_\nu^{(1)})(z)=Ce_{q^2}(-(1-q^2)^2z^2).
$$
By setting $z=0$, we have
$$
C=\frac{q^{-\nu^2}(1-q^2)}\pi
$$
and
$$
W(J_\nu^{(1)},Y_\nu^{(1)})(z)=\frac{q^{-\nu^2}(1-q^2)}\pi
e_{q^2}(-(1-q^2)^2z^2).
$$

The proof of the Proposition for $\de=1 ~{\rm and} ~0$ is similar. $\Box$

\begin{predl}\label{p2.2}
The functions $I_\nu^{(j)}(2(1-q^2)z;q^2)$ and $K_\nu^{(j)}(2(1-q^2)z;q^2)$
satisfy the difference equation
$$
f(q^{-1}z)-(q^{-\nu}+q^\nu)f(z)+f(qz)=q^{-\de}(1-q^2)^2z^2f(q^{1-\de}z),
$$
and form the fundamental system of solutions of this equation with
$q^2$-Wronskian
$$
W(I_\nu^{(j)},K_\nu^{(j)})=\left\{\begin{array}{lcl}
\frac{q^{-\nu^2}(1-q^2)}2 e_{q^2}((1-q^2)^2z^2)&for&\de=2\\
 \frac{q^{-\nu^2}(1-q^2)}2&for&\de=1\\
\frac{q^{-\nu^2}(1-q^2)}2 E_{q^2}(-(1-q^2)^2z^2)&for&\de=0.\\
\end{array}\right.
$$
\end{predl}
The {\it Proof} is the same as above.

\section{The Laurent series associated with the $q$-exponentials}
\setcounter{equation}{0}

Let $z$ be a real positive number, and $u=e^{i\te}z(1-q^2)$. Consider the
products
\beq{3.1}
\La^{(j)}(u)=e_q^{(j)}(u)e_q^{(j)}(\frac qu); ~~~~~j=1, ~2, ~3.
\eq

\begin{predl}\label{p3.1}
The functions $\La^{(j)}(u)$ (\ref{3.1}) can be represented as the Laurent
series
\beq{3.2}
\La^{(j)}(u)=\sum_{l=-\infty}^\infty u^lq^{\frac{2-\de}4l^2-\frac l2}
I_l^{(j)}(2q^{\frac\de4};q).
\eq
If $j=1, ~q<|u|<1.$
\end{predl}
{\it Proof.}  Let $j=1$ and $q<|u|<1.$ Substitute (\ref{1.2}) into
(\ref{3.1}).
$$
\La^{(1)}(u)=\sum_{n=0}^\infty\frac{u^n}{(q,q)_n}
\sum_{m=0}^\infty\frac{q^mu^{-m}}{(q,q)_m}=
\sum_{l=1}^\infty u^{-l}\sum_{k=0}^\infty\frac{q^{k+l}}{(q,q)_k(q,q)_{k+l}}+
\sum_{l=0}^\infty u^l\sum_{k=0}^\infty\frac{q^k}{(q,q)_k(q,q)_{k+l}}=
$$
\beq{3.3}
\sum_{l=1}^\infty\frac{q^lu^{-l}}{(q,q)_l}
\sum_{k=0}^\infty\frac{q^k}{(q,q)_k(q^{l+1},q)_k}+
\sum_{l=0}^\infty\frac{u^l}{(q,q)_l}
\sum_{k=0}^\infty\frac{q^k}{(q,q)_k(q^{l+1},q)_k}.
\eq
It follows from (\ref{2.1}) that if $q<|u|<1$
$$
\sum_{k=0}^\infty\frac{q^k}{(q,q)_k(q^{l+1},q)_k}=
(q,q)_lq^{-\frac{l}{2}}I_l^{(1)}(2\sqrt q;q).
$$
Since $I_{-l}^{(j)}(z;q)=I_l^{(j)}(z;q)$ we have
$$
\La^{(1)}(u)=\sum_{l=-\infty}^\infty u^lq^{-\frac l2}I_l^{(1)}(2\sqrt q;q).
$$

If $j=2$ or $3$ for any $u\ne 0$ we have
$$
\La^{(j)}(u)=\sum_{n=0}^\infty\frac{q^{\frac{2-\de}4n(n-1)}u^n}{(q,q)_n}
\sum_{m=0}^\infty\frac{q^{\frac{2-\de}4m(m-1)+m}u^{-m}}{(q,q)_m}=
$$
$$
\sum_{l=1}^\infty u^{-l}
\sum_{k=0}^\infty\frac{q^{\frac{2-\de}4k(k-1)+\frac{2-\de}4(k+l)(k+l-1)+k+l}}
{(q,q)_k(q,q)_{k+l}}+
\sum_{l=0}^\infty u^l
\sum_{k=0}^\infty\frac{q^{\frac{2-\de}4k(k-1)+\frac{2-\de}4(k+l)(k+l-1)+k}}
{(q,q)_k(q,q)_{k+l}}=
$$
\beq{3.4}
\sum_{l=1}^\infty\frac{q^{\frac{2-\de}4l(l-1)+l}u^{-l}}{(q,q)_l}
\sum_{k=0}^\infty\frac{q^{\frac{2-\de}2k(k+l)+\frac\de2k}}{(q,q)_k(q^{l+1},q)_k}+
\sum_{l=1}^\infty\frac{q^{\frac{2-\de}4l(l-1)}u^l}{(q,q)_l}
\sum_{k=0}^\infty\frac{q^{\frac{2-\de}2k(k+l)+\frac\de2k}}{(q,q)_k(q^{l+1},q)_k}.
\eq
It follows from (\ref{2.1}) that for any $u$
$$
\sum_{k=0}^\infty\frac{q^{\frac{2-\de}2k(k+l)+k\frac\de2}}{(q,q)_k(q^{l+1},q)_k}=
(q,q)_lq^{-l\frac\de4}I_l^{(j)}(2q^{\frac\de4};q).
$$
So
\beq{3.5}
\La^{(j)}(u)=\sum_{l=-\infty}^\infty u^lq^{\frac{2-\de}4l^2-\frac l2}
I_l^{(j)}(2q^{\frac\de4};q).
\eq

Obviously for $j=1 ~(\de=2)$ we get (\ref{3.2}) from (\ref{3.5}). $\Box$

\begin{cor}\label{c3.1}  If $u$ is real and $q<u<1$
\beq{3.6}
\La^{(2)}(u)\le\La^{(3)}(u)\le\La^{(1)},
\eq
and for $u>1$
\beq{3.7}
\La^{(2)}(u)\le\La^{(3)}(u).
\eq
\end{cor}
The {\it Proof}  follows from the comparison of the coefficients in
(\ref{3.3}) and (\ref{3.4}) for $j=1, ~2, ~3 ~~(\de=2, ~0, ~1)$.  $\Box$

\section{The behaviour of the $q$-exponentials for the large value of the argument}
\setcounter{equation}{0}

\begin{predl}\label{p4,1}
The function $\La^{(1)}(u)$ (\ref{3.1}) for $j=1$ satisfies the equation
\beq{4.1}
\La^{(1)}(qu)+u\La^{(1)}(u)=0
\eq
and has the form
\beq{4.2}
\La^{(1)}(u)=u^{\frac12\log_qu-\frac12+\frac{i\pi}{\ln q}}.
\eq
\end{predl}
{\it Proof.}  Consider the function
$\La^{(1)}=\frac1{(u,q)_\infty(\frac qu,q)_\infty}$.
$$
D_z\La^{(1)}(u)e^{-i\te}=u^{-1}\left(\La^{(1)}(u)-\La^{(1)}(qu)\right)=
\frac1u\left(\frac1{(u,q)_\infty(\frac qu,q)_\infty}-
\frac1{(qu,q)_\infty(\frac1u,q)_\infty}\right)=
$$
$$
\frac{(1-\frac1u)}{u(1-u)}\frac1{(qu,q)_\infty(\frac1u,q)_\infty}-
\frac{1-u}{u(1-\frac1u)}\frac1{(u,q)_\infty(\frac qu,q)_\infty}=
-u^{-2}\La^{(1)}(qu)+\La^{(1)}(u)
$$

The last equality gives us (\ref{4.1}).

The fact that (\ref{4.2}) satisfies (\ref{4.1}) can be checked directly. $\Box$

\begin{predl}\label{p4,2}
The function $\La^{(2)}(u)$ (\ref{3.1}) for $j=2$ satisfies the equation
\beq{4.3}
u\La^{(2)}(qu)-\La^{(2)}(u)=0
\eq
and has the form
\beq{4.4}
\La^{(2)}(u)=u^{-\frac12\log_qu+\frac12}.
\eq
\end{predl}
{\it Proof.}  Consider the function
$\La^{(2)}=(-u,q)_\infty(-\frac qu,q)_\infty$.
$$
D_z\La^{(2)}(u)e^{-i\te}=u^{-1}\left(\La^{(2)}(u)-\La^{(2)}(qu)\right)=
\frac1u\left((-u,q)_\infty(-\frac qu,q)_\infty-
(-qu,q)_\infty(-\frac1u,q)_\infty\right)=
$$
$$
\frac{1+u}{u(1+\frac1u)}((-qu,q)_\infty(-\frac1u,q)\infty-
\frac{1+\frac1u}{u(1+u)}(-u,q)_\infty(-\frac qu,q)\infty)=
$$
$$
\La^{(2)}(qu)-u^{-2}\La^{(2)}(u).
$$
The last equality gives us (\ref{4.3})

The fact that (\ref{4.4}) satisfies (\ref{4.3}) can be checked directly. $\Box$

\bigskip
Represent the function $\La^{(3)}(u)$ (\ref{3.1}) for $j=3$ by product
\beq{4.5}
\La^{(3)}(u)= \La(u)\La(\frac qu)
\eq
and let
$$
D_z\La(u)=\La(q^{\frac12}u),
~~~~~D_z\La(\frac qu)=-\frac1{u^2}\La(\frac{q^{\frac12}}u).
$$
Then
$$
D_z(\La(u)\La(\frac qu))=-\frac1{u^2}\La(u)\La(\frac{q^{\frac12}}u)+
\La(q^{\frac12}u)\La(\frac1u).
$$

\begin{predl}\label{p4,3}
The function $\La(u)$ satisfies the equation
\beq{4.6}
\La(u)\La(\frac qu)-\La(qu)\La(\frac1u)-u\La(q^{\frac12}u)\La(\frac1u)+
\frac1u\La(u)\La(\frac{q^{\frac12}}u)=0
\eq
and has the forms
\beq{4.7}
\La_1(u)=u^{\log_qu-\frac32+\frac{2i\pi}{\ln q}}
\eq
or
\beq{4.8}
\La_2(u)=u^{-\frac13\log_qu+\frac12}.
\eq
\end{predl}
{\it Proof.}  Obviously
$$
D_z\La^{(3)}(u)e^{-i\te}=u^{-1}\left[\La(u)\La(\frac qu)-\La(qu)\La(\frac1u)\right]=
\La(q^{\frac12}u)\La(\frac1u)-u^{-2}\La(u)\La(\frac{q^{\frac12}}u).
$$
So we have (\ref{4.6}).

The fact that (\ref{4.7}) and (\ref{4.8}) satisfy (\ref{4.6}) can be checked
directly. $\Box$

\bigskip
Now investigate the behaviour of the $q$-exponentials for the large value of the
argument. Consider the cases $j=1, ~2, ~3 ~(\de=2, ~0, ~1)$ separately.

1. $j=1.$ Consider (\ref{4.2}). Put $u=e^{i\te}|u|, ~|u|=q^{n+\la},
~n=\left[\frac{\ln|u|}{\ln q}\right], ~[a]$ - is an integral
part of $a$, and $0\le\la<1.$ Then
$$
\La^{(1)}(u)=C^{(1)}e^{i(\te+\pi)n}q^{\frac12n(n-1)+\la n},
$$
where $C^{(1)}$ is an independent on $n$ constant. Assuming $n=0$
we obtain
$$
C^{(1)}=\La^{(1)}(e^{i\te}q^\la).
$$

On the other hand it follows from (\ref{1.2}) that
$$
\La^{(1)}(u)=\frac1{(u,q)_\infty(qu^{-1},q)_\infty}=
\frac{1}{(q^{n+\la}e^{i\te},q)_\infty
(q^{1-n-\la}e^{-i\te},q)_\infty}=
$$
$$
\frac{(q^\la e^{i\te},q)_n}
{(q^\la e^{i\te},q)_\infty(q^{1-n-\la}e^{-i\te},q)_n
(q^{1-\la}e^{-i\te},q)_\infty}=
\frac{q^{\la n+\frac12n(n-1)}e^{in\la(\te+\pi)}}
{(q^\la e^{i\te},q)_\infty(q^{1-\la}e^{-i\te},q)_\infty}.
$$
Thus
$$
C^{(1)}=\frac{1}{(q^\la e^{i\te},q)_\infty(q^{1-\la}e^{-i\te},q)_\infty}
$$
and for any $u\ne 0$ \cite{I}
\beq{4.9}
e_q^{(1)}(u)=\frac{\La^{(1)}(u)}{e_q^{(1)}(\frac qu)}=
\frac{q^{\la n)+\frac12n(n-1)}e^{i(\pi+\te)n}
(q^{1-n-\la}e^{-i\te},q)_\infty}
{(q^\la e^{i\te},q)_\infty(q^{1-\la}e^{-i\te},q)_\infty}.
\eq

2. $j=2.$ Consider (\ref{4.4})
$$
\La^{(2)}(u)=u^{-\frac12(\log_qu-1)}
$$
and assume $|u|=q^{n+\la}$ again. Then
$$
\La^{(2)}(u)=C^{(2)}
e^{-i\te n}q^{-\frac12n(n-1)-\la n)}.
$$
On the other hand it follows from (\ref{1.3}) that
$$
\La^{(2)}(u)=(-u,q)_\infty(-qu^{-1},q)_\infty=
(-q^{n+\la}e^{i\te},q)_\infty(-q^{1-n-\la}e^{-i\te},q)_\infty=
$$
$$
\frac{(-q^\la e^{i\te},q)_\infty(-q^{1-n-\la}e^{-i\te},q)_n
(-q^{1-\la}e^{-i\te},q)_\infty}{(-q^\la e^{i\te},q)_n}=
$$
$$
q^{-\la n-\frac12n(n-1)}e^{-i\te n}
(-q^\la e^{i\te},q)_\infty(-q^{1-\la}e^{-i\te},q)_\infty,
$$
Thus
$$
C^{(2)}=(-q^\la e^{i\te},q)_\infty(-q^{1-\la}e^{-i\te},q)_\infty
$$
and for any $u\ne 0$
\beq{4.10}
e_q^{(2)}(u)=\frac{\La^{(2)}(u)}{e_q^{(2)}(\frac q{u})}=
q^{-\la n-\frac12n(n-1)}e^{-i\te n}
\frac{(-q^\la e^{i\te},q)_\infty(-q^{1-\la}e^{-i\te},q)_\infty}
{(-q^{1-n-\la}e^{-i\te},q)_\infty}.
\eq

3. $j=3.$ The equation (\ref{4.6}) has two solutions (\ref{4.7}) and
(\ref{4.8}). Substituting these functions into (\ref{4.5}) we obtain
$$
\La_1^{(3)}(u)=C_1u^{2(\log_qu-1)},
~~~~~~~~~\La_2^{(3)}(u)=C_2u^{-\frac23(\log_qu-1)}.
$$

For an arbitrary real $u>1$ take the solution $\La_2^{(3)}(u)$ and require
that
$$
u^{-\frac12(\log_qu-1)}\le q^\mu u^{-\frac23(\log_qu-1)},
$$
according to (\ref{3.7}).
$$
\mu-\frac23(\log_qu-1)\log_qu\le-\frac12(\log_qu-1)\log_qu,
$$
$$
\mu\le\frac16(\log_qu-1)\log_qu,
$$
and we put $\mu=-\frac1{24}$.  Then (\ref{3.6}) is fulfilled for $u: ~q<u<1.$
So for an arbitrary complex $u\ne 0$
$$
\La^{(3)}(u)=q^{-\frac1{24}}u^{-\frac23(\log_qu-1)},
$$
and we assume $|u|=q^{n+\la}$ again.
$$
\La^{(3)}(u)=C^{(3)}q^{-\frac23n(n-1)-\frac43n\la}
e^{-\frac{4i\te}3n},
$$
$$
C^{(3)}=q^{-\frac1{24}}e_q^{(3)}(q^{\{\la\}}e^{i\te})
e_q^{(3)}(q^{1-\{\la\}}e^{-i\te}).
$$
Hence for any $u\ne 0$
\beq{4.11}
e_q^{(3)}(u)=\frac{\La^{(3)}(u)}{e_q^{(3)}(\frac{q}u)}=
q^{-\frac23n(n-1)-\frac43n\la-\frac1{24}}e^{-\frac{4i\te}3n}
\frac{e_q^{(3)}(q^\la e^{i\te})e_q^{(3)}(q^{1-\la}e^{-i\te})}
{e_q^{(3)}(q^{1-n-\la}e^{-i\te})}.
\eq

\begin{predl}\label{p4.4}
The asymptotic behaviour of the $q$-exponentials has the form
$$
e_q^{(j)}(u)=\left\{\begin{array}{ccc}
q^{\frac12N}e^{i(\te+\pi)n}C^{(1)}(1+o(q^{-n})) &{\rm for} &j=1\\
q^{-\frac12N}e^{-i\te n}C^{(2)}(1+o(q^{-n}))& {\rm for} &j=2\\
q^{-\frac23N-\frac1{24}}
e^{-\frac{4i\te}3n}C^{(3)}(1+o(q^{-n}))& {\rm for} &j=3,\\
\end{array}\right.
$$
where
\beq{4.12}
N=n(n-1)+2\la n,
\eq
and
$$
C^{(j)}=e_q^{(j)}(q^\la e^{i\te})e_q^{(j)}(q^{1-\la}e^{-i\te}).
$$
\end{predl}
{\it Proof.} Obviously
$$
\lim_{n\to-\infty}e_q^{(j)}(q^{1-n-\la}e^{-i\te})=1.
$$
Now the statement of the Proposition follows from (\ref{4.9}) -
(\ref{4.11}).
$\Box$

\section{The representations of the $q^2$-Bessel functions as the Laurent series}
\setcounter{equation}{0}

Let $u$ be real.

In \cite{OR1} it has been received the following representations for the
modified $q^2$-Bessel functions of kind 1 and 2
\beq{5.1}
I_\nu^{(j)}(2u;q^2)=\frac{a_\nu}{\sqrt{2u}}
\left[e_q^{(j)}(u)\Phi_\nu(u)+ie^{i\nu\pi}e_q^{(j)}(-u)\Phi_\nu(-u)\right],
\eq
\beq{5.2}
K_\nu^{(j)}(2u;q^2)=\frac{q^{-\nu^2+\frac12}(1-q^2)}
{2a_\nu\sqrt{2u}}e_q^{(j)}(-u)\Phi_\nu(-u), ~~{\rm  for} ~\nu\ne n,
\eq
where
$$
u=(1-q^2)z, ~~~~a_\nu^2=a_{-\nu}^2=\left\{\begin{array}{lcl}
\frac{q^{-\nu+\frac12}(1-q^2)}
{2\G_{q^2}(\nu)\G_{q^2}(1-\nu)\sin \nu\pi}&{\rm for}&\nu\ne n\\
\frac{q^{-n^2+\frac12}\ln q^{-2}}{2\pi}&{\rm for}&\nu=n
\end{array}\right.,
$$
\beq{5.3}
\Phi_\nu(u)=
\phantom._2\Phi_1(q^{\nu+\frac12},q^{-\nu+\frac12};-q;q,\frac q{u}).
\eq
In the case $\nu=n ~~K_n^{(j)}$ can be received as
the limit of $K_\nu^{(j)}$ if $\nu$ tends to $n$.

It follows from (\ref{1.2}), (\ref{1.3}) and (\ref{5.3}) that
$\sqrt uI_\nu^{(1)}(2u;q^2)$ is a meromorphic function and has simple poles
at the points $u=\pm q^{-n}, ~~n=0, 1,...$.
$\sqrt uK_\nu^{(1)}(2u;q^2)$ is a meromorphic function and has simple poles
at the points $u=-q^{-n}, ~~n=0, 1,...$.
$~~\sqrt uI_\nu^{(2)}(2u;q^2)$ and $~~\sqrt uK_\nu^{(2)}(2u;q^2)$ are the
holomorphic functions outside a neighborhood of zero (see \cite{OR1})

Using (\ref{2.2}) we can write
\beq{5.4}
J_\nu^{(j)}(2u;q^2)=\frac{a_\nu}{\sqrt{2u}}
\left[e^{-i(\frac\pi4+\nu\frac\pi2)}e_q^{(j)}(iu)\Phi_\nu(iu)+
ie^{i(\frac\pi4+\nu\frac\pi2)}e_q^{(j)}(-iu)\Phi_\nu(-iu)\right],
\eq
\beq{5.5}
Y_\nu^{(j)}(2u;q^2)=
\eq
$$
-\frac{q^{-\nu^2+\frac12}(1-q^2)}{\pi\sqrt{2u}a_\nu}
\left[e^{i(\frac\pi4-\nu\frac\pi2)}e_q^{(j)}(iu)\Phi_\nu(iu)+
ie^{-i(\frac\pi4-\nu\frac\pi2)}e_q^{(j)}(-iu)\Phi_\nu(-iu)\right],
$$
if $\nu$ is not integer. In the case $\nu=n ~~~Y_n^{(j)}$ can be received as
the limit of $Y_\nu^{(j)}$ if $\nu$ tends to $n$.

The functions $\sqrt uJ_\nu^{(1)}(2u;q^2)$ and $\sqrt uY_\nu^{(1)}(2u;q^2)$
are the meromorphic functions and have simple poles at the points
$u=\pm iq^{-n}, ~~n=0, 1,...$. The functions $\sqrt uJ_\nu^{(2)}(2u;q^2)$
and $\sqrt uY_\nu^{(2)}(2u;q^2)$ are the holomorphic functions outside a
neighborhood of zero.

To receive the similar formulas for $q^2$-Bessel functions of kind 3 we
use the following

\begin{lem}\label{l5.1}
Let the functions
\beq{5.6}
f^{(j)}(z)=\sum_{k=-\infty}^\infty b_k^{(j)}z^{k-\frac12}
\eq
satisfy the equations
\beq{5.7}
q^{\frac12}f^{(j)}(q^{-1}z)-(q^{-\nu}+q^\nu)f^{(j)}(z)+q^{-\frac12}f^{(j)}(qz)=
q^{-\frac{\de+1}2}z^2f^{(j)}(q^{1-\de}z)
\eq
where $j=-\frac32\de^2+\frac52\de+2, ~~~\de=0, 1$ and $2$. Then
\beq{5.8}
b_k^{(3)}=\sqrt{b_k^{(1)}b_k^{(2)}}
\eq
for any $k$
\end{lem}

{\it Proof}. It follows from (\ref{5.6}) and (\ref{5.7}) that
$$
b_k^{(j)}=\frac{b_{k-2}^{(j)}q^{k\de-3(1-\frac\de2)}}
{(1-q^{-\nu+k-\frac12})(1-q^{\nu+k-\frac12})}=
\left\{\begin{array}{ccc}
\frac{b_0^{(j)}q^{l(l\de-3+\frac52\de)}}
{(q^{-\nu+\frac32},q^2)_l(q^{\nu+\frac32},q^2)_l}& {\rm for}& k=2l\\
\frac{b_1^{(j)}q^{l(l\de-3+2\de)}}
{(q^{-\nu+\frac32},q^2)_l(q^{\nu+\frac32},q^2)_l}& {\rm for}& k=2l+1,\\
\end{array}\right.
$$
i.e. (\ref{5.8}) is fulfilled. $\Box$

Find the coefficient of function
$$
\frac1{\sqrt u}e_q^{(j)}(u)\Phi_\nu(u)=\frac1{\sqrt u}
\left[\sum_{l=1}^\infty c_{l-}^{(j)}u^{-l}+
\sum_{l=0}^\infty c_{l+}^{(j)}u\right], ~~~~~j=1, 2,
$$
in their decompositions in the Laurent series. Because these functions satisfy
(\ref{5.7}) for $\de=2$ and $0$ respectively the following series satisfies
(\ref{5.7}) for $\de=1$
\beq{5.9}
\left[\sum_{l=1}^\infty\sqrt{c_{l-}^{(1)}c_{l-}^{(2)}}u^{-l-\frac12}+
\sum_{l=0}^\infty\sqrt{c_{l+}^{(1)}c_{l+}^{(2)}}u^{l-\frac12}\right].
\eq
Using (\ref{1.2}), (\ref{1.3}) and (\ref{5.3}), we can write for $j=1, ~2$
$$
c_{l-}^{(j)}=\sum_{k=0}^\infty\frac{q^{\frac{2-\de}4k(k-1)}}{(q,q)_k}
\frac{(q^{-\nu+\frac12},q)_{k+l}(q^{\nu+\frac12},q)_{k+l}}
{(q^2,q^2)_{k+l}}q^{k+l}=
$$
$$
\frac{(q^{-\nu+\frac12},q)_l(q^{\nu+\frac12},q)_l}{(q^2,q^2)_l}q^l
\sum_{k=0}^\infty\frac{(q^{-\nu+l+\frac12},q)_k(q^{\nu+l+\frac12},q)_k}
{(q^{2l+2},q^2)_k(q,q)_k}q^{\frac{2-\de}4k(k-1)+k}.
$$
$$
c_{l+}^{(j)}=\sum_{k=0}^\infty\frac{q^{\frac{2-\de}4(k+l)(k+l-1)}}{(q,q)_{k+l}}
\frac{(q^{-\nu+\frac12},q)_k(q^{\nu+\frac12},q)_k}
{(q^2,q^2)_k}q^k=
$$
$$
\frac{q^{\frac{l(l-1)}2}}{(q,q)_l}
\sum_{k=0}^\infty\frac{(q^{-\nu+\frac12},q)_k(q^{\nu+\frac12},q)_k}
{(q^2,q^2)_k(q^{l+1},q)_k}q^{\frac{2-\de}4k(k+l-1)+k}.
$$
For this reason we have
\beq{5.10}
\sqrt{c_{l-}^{(1)}c_{l-}^{(2)}}=\frac{(q^{-\nu+\frac12},q)_l(q^{\nu+\frac12},q)_l}
{(q^2,q^2)_l}q^l\times
\eq
$$
\times\sqrt{\sum_{k=0}^\infty\frac{(q^{-\nu+l+\frac12},q)_k(q^{\nu+l+\frac12},q)_k}
{(q^{2l+2},q^2)_k(q,q)_k}q^k
\sum_{m=0}^\infty\frac{(q^{-\nu+l+\frac12},q)_m(q^{\nu+l+\frac12},q)_m}
{(q^{2l+2},q^2)_m(q,q)_m}q^{\frac{m(m+1)}2}},
$$
\beq{5.11}
\sqrt{c_{l+}^{(1)}c_{l+}^{(2)}}=
\eq
$$
\frac{q^{\frac{l(l-1)}4}}{(q,q)_l}
\sqrt{\sum_{k=0}^\infty\frac{(q^{-\nu+\frac12},q)_k(q^{\nu+\frac12},q)_k}
{(q^2,q^2)_k(q^{l+1},q)_k}q^k
\sum_{m=0}^\infty\frac{(q^{-\nu+\frac12},q)_m(q^{\nu+\frac12},q)_m}
{(q^2,q^2)_m(q^{l+1},q)_m}q^{\frac{m(m+l+1)}2}}.
$$

The representations of $I_\nu^{(j)}, ~~j=1, ~2,$ have the form
$$
I_\nu^{(j)}(2u;q^2)=\frac{a_\nu}{\sqrt{2u}}\times
$$
$$
\times\left[\sum_{l=1}^\infty c_{l-}^{(j)}u^{-l}+
\sum_{l=0}^\infty c_{l+}^{(j)}u^l+
ie^{i\nu\pi}\left(\sum_{l=1}^\infty (-1)^lc_{l-}^{(j)}u^{-l}+
\sum_{l=0}^\infty (-1)^lc_{l+}^{(j)}u^l\right)\right]=
$$
$$
\frac{a_\nu}{\sqrt{2u}}\left[\sum_{l=1}^\infty(1+ie^{i(\nu+l)\pi})
c_{l-}^{(j)}u^{-l}+\sum_{l=0}^\infty(1+ie^{i(\nu+l)\pi})
c_{l+}^{(j)}u^l\right].
$$
Hence
\beq{5.12}
I_\nu^{(3)}(2u;q^2)=\frac{a_\nu}{\sqrt{2u}}\times
\eq
$$
\times\left[\sum_{l=1}^\infty(1+ie^{i(\nu+l)\pi})
\sqrt{c_{l-}^{(1)}c_{l-}^{(2)}}u^{-l}
+\sum_{l=0}^\infty(1+ie^{i(\nu+l)\pi})
\sqrt{c_{l+}^{(1)}c_{l+}^{(2)}}u^l\right].
$$
In the similar way we obtain
\beq{5.13}
K_\nu^{(3)}(2u;q^2)=\frac{q^{-\nu^2+\frac12}(1-q^2)}{2a_\nu\sqrt{2u}}
\left[\sum_{l=1}^\infty(-1)^l\sqrt{c_{l-}^{(1)}c_{l-}^{(2)}}u^{-l}
+\sum_{l=0}^\infty(-1)^l\sqrt{c_{l+}^{(1)}c_{l+}^{(2)}}u^l\right],
\eq
\beq{5.14}
J_\nu^{(3)}(2u;q^2)=\frac{a_\nu}{\sqrt{2u}}e^{-i\frac\pi4-i\frac\nu2\pi}\times
\eq
$$
\times\left[\sum_{l=1}^\infty(1+ie^{i\nu\pi})
\sqrt{c_{l-}^{(1)}c_{l-}^{(2)}}u^{-l}
+\sum_{l=0}^\infty(1+ie^{i\nu\pi})
\sqrt{c_{l+}^{(1)}c_{l+}^{(2)}}u^l\right],
$$
\beq{5.15}
Y_\nu^{(3)}(2u;q^2)=\frac{q^{-\nu^2+\frac12}(1-q^2)}{a_\nu\pi\sqrt{2u}}
e^{i\frac\pi4-i\frac\nu2\pi}\times
\eq
$$
\times\left[\sum_{l=1}^\infty(-1+ie^{i\nu\pi})
\sqrt{c_{l-}^{(1)}c_{l-}^{(2)}}u^{-l}+\sum_{l=0}^\infty(-1+ie^{i\nu\pi})
\sqrt{c_{l+}^{(1)}c_{l+}^{(2)}}u^l\right].
$$

\section{The asymptotic behaviour of the $q^2$-Bessel functions of type 1 and 2}
\setcounter{equation}{0}

Let $u=q^{n+\la}, ~~n=\left[\frac{\ln u}{\ln q}\right],$
where $[a]$ is an integral part of $a$ and $0\le\la<0$.
Introduce the designations
\beq{6.1}
Q_\pm^{(j)}=\frac{e_q^{(j)}(\pm q^\la)
e_q^{(j)}(\pm q^{1-\la})}{e_q^{(j)}(\pm q^{1-n-\la})},
~~~Q_{\pm i}^{(j)}=\frac{e_q^{(j)}(\pm iq^\la)
e_q^{(j)}(\mp iq^{1-\la})}{e_q^{(j)}(\mp iq^{1-n-\la})}
\eq
and define $N$ by (\ref{4.12}) as previously. Then it follows from
(\ref{4.9}) - (\ref{4.11}) that
$$
e_q^{(1)}(\pm u)=q^{\frac12N}e^{i\frac\pi2(1\pm1)n}Q_\pm^{(1)},
~~~~e_q^{(1)}(\pm iu)=q^{\frac12N}e^{i\pi(1\pm\frac12)n}Q_{\pm i}^{(1)}.
$$
$$
e_q^{(2)}(\pm u)=q^{-\frac12N}e^{i\frac\pi2(1\mp1)n}Q_\pm^{(2)},
~~~~e_q^{(2)}(\pm iu)=q^{-\frac12N}e^{\mp i\frac\pi2n}Q_{\pm i}^{(2)}.
$$
Using (\ref{5.1}) and (\ref{5.2}) we obtain
\beq{6.2}
I_\nu^{(1)}(2u;q^2)=\frac{a_\nu}{\sqrt{2u}}
q^{\frac12N}\left[e^{in\pi}Q_+^{(1)}\Phi_\nu(u)+ie^{i\nu\pi}
Q_-^{(1)}\Phi_\nu(-u)\right],
\eq
\beq{6.3}
K_\nu^{(1)}(2u;q^2)=\frac{q^{-\nu^2+\frac12}(1-q^2)}{2a_\nu\sqrt{2u}}
q^{\frac12N}Q_-^{(1)}\Phi_\nu(-u),
\eq
\beq{6.4}
I_\nu^{(2)}(2u;q^2)=\frac{a_\nu}{\sqrt{2u}}
q^{-\frac12N}\left[Q_+^{(2)}\Phi_\nu(u)+
ie^{i\nu\pi+in)\pi}Q_-^{(2)}\Phi_\nu(-u)\right],
\eq
\beq{6.5}
K_\nu^{(2)}(2u);q^2)=\frac{q^{-\nu^2+\frac12}(1-q^2)}{2a_\nu\sqrt{2u}}
q^{-\frac12N}e^{i\pi n} Q_-^{(2)}\Phi_\nu(-u).
\eq

Us the formulas (\ref{5.4}), (\ref{5.5}) and (\ref{6.1}). Then
\beq{6.6}
J_\nu^{(1)}(2u;q^2)=\frac{a_\nu}{\sqrt{2u}}q^{\frac12N}\times
\eq
$$
\times\left[e^{-i(2\nu+1)\frac\pi4+i\frac{3\pi}2n}
Q_{+i}^{(1)}\Phi_\nu(iu)+e^{i(2\nu+1)\frac\pi4+i\frac\pi2n}
Q_{-i}^{(1)}\Phi_\nu(-iu)\right]
$$
\beq{6.7}
Y_\nu^{(1)}(2u;q^2)=
\frac{q^{-\nu^2+\frac12}(1-q^2)}{\pi a_\nu\sqrt{2u}}q^{\frac12N}\times
\eq
$$
\times\left[-e^{-i(2\nu-1)\frac\pi4+i\frac{3\pi}2n}
Q_{+i}^{(1)}\Phi_\nu(iu)-e^{i(2\nu-1)\frac\pi4+i\frac\pi2n}
Q_{-i}^{(1)}\Phi_\nu(-iu)\right]
$$
\beq{6.8}
J_\nu^{(2)}(2u;q^2)=\frac{a_\nu}{\sqrt{2u}}q^{-\frac12N}\times
\eq
$$
\times\left[e^{-i(2\nu+1)\frac\pi4-i\frac\pi2n}
Q_{+i}^{(2)}\Phi_\nu(iu)+e^{i(2\nu+1)\frac\pi4+i\frac\pi2n}
Q_{-i}^{(2)}\Phi_\nu(-iu)\right]
$$
\beq{6.9}
Y_\nu^{(2)}(2u;q^2)=\frac{q^{-\nu^2+\frac12}(1-q^2)}
{\pi a_\nu\sqrt{2u}}q^{-\frac12N}\times
\eq
$$
\times\left[-e^{-i(2\nu-1)\frac\pi4-i\frac\pi2n}
Q_{+i}^{(2)}\Phi_\nu(iu)-e^{i(2\nu-1)\frac\pi4+i\frac\pi2n}
Q_{-i}^{(2)}\Phi_\nu(-iu)\right].
$$

Now we can describe the asymptotic behaviour of $q^2$-Bessel functions
of type $1$ and $2$. Note that if $n\to-\infty ~~(|u|\to\infty)$,
$$
\Phi_\nu(\pm u)Q_{\pm}^{(j)}\to e_q^{(j)}(\pm q^\la)
e_q^{(j)}(\pm q^{1-\la})
$$
and
$$
\Phi_\nu(\pm iu)Q_{\pm i}^{(j)}\to\frac{1}e_q^{(j)}(\pm iq^\la)
e_q^{(j)}(\mp iq^{1-\la}).
$$

\begin{predl}\label{p6.1}  The following asymptotic formulas take place for
$I_\nu^{(j)}(2u;q^2), ~K_\nu^{(j)}(2u;q^2), \\
J_\nu^{(1)}(2u;q^2), ~Y_\nu^{(1)}(2u;q^2), ~~j=1, ~2$
$$
I_\nu^{(1)}(2u;q^2)=\frac{a_\nu}{\sqrt{2u}}q^{\frac12N}\times
$$
$$
\times\left[[e^{in\pi}e_q^{(1)}(q^\la)e_q^{(1)}(q^{1-\la})+
ie^{i\nu\pi}e_q^{(1)}(-q^\la)e_q^{(1)}(-q^{1-\la})\right](1+o(q^{-n})),
$$
$$
K_\nu^{(1)}(2u;q^2)=\frac{q^{-\nu^2+\frac12}(1-q^2)}
{2a_\nu\sqrt{2u}}q^{\frac12N}e_q^{(1)}(-q^\la)
e_q^{(1)}(-q^{1-\la})(1+o(q^{-n})),
$$
$$
I_\nu^{(2)}(2u;q^2)=\frac{a_\nu}{\sqrt{2u}}q^{-\frac12N}
\left[e_q^{(2)}(q^\la)e_q^{(2)}(q^{1-\la})+\right.
$$
$$
\left.ie^{i\nu\pi+i\pi n}e_q^{(2)}(-q^\la)
e_q^{(2)}(-q^{1-\la})\right](1+o(q^{-n})),
$$
$$
K_\nu^{(2)}(2u;q^2)=
\frac{q^{-\nu^2+\frac12}(1-q^2)}{2a_\nu\sqrt{2u}}q^{-\frac12N}
e^{i\pi n}e_q^{(2)}(-q^\la)e_q^{(2)}(-q^{1-\la})(1+o(q^{-n})),
$$
$$
J_\nu^{(1)}(2u;q^2)=\frac{a_\nu}{\sqrt{2u}}q^{\frac12N}
\left[e^{-i(2\nu+1)\frac\pi4+i\frac{3\pi}2n}
e_q^{(1)}(iq^\la)e_q^{(1)}(-iq^{1-\la})+\right.
$$
$$
\left.e^{i(2\nu+1)\frac\pi4+i\frac\pi2(n+\la)}
e_q^{(1)}(-iq^\la)e_q^{(1)}(iq^{1-\la})\right](1+o(q^{-n})),
$$
$$
Y_\nu^{(1)}(2u;q^2)=\frac{q^{-\nu^2+\frac12}(1-q^2)}{\pi\sqrt{2u}}q^{\frac12N}
\left[-e^{-i(2\nu-1)\frac\pi4+i\frac{3\pi}2n}
e_q^{(1)}(iq^\la)e_q^{(1)}(-iq^{1-\la})-\right.
$$
$$
\left.e^{i(2\nu-1)\frac\pi4+i\frac\pi2n}
e_q^{(1)}(-iq^\la e_q^{(1)}(iq^{1-\la})\right](1+o(q^{-n})),
$$
$$
J_\nu^{(2)}(2uz;q^2)=\frac{a_\nu}{\sqrt{2u}}q^{-\frac12N}
\left[e^{-i(2\nu+1)\frac\pi4-i\frac{\pi}2n}
e_q^{(2)}(iq^\la)e_q^{(2)}(-iq^{1-\la})+\right.
$$
$$
\left.e^{i(2\nu+1)\frac\pi4+i\frac\pi2n}
e_q^{(2)}(-iq^\la)e_q^{(2)}(iq^{1-\la})\right](1+o(q^{-n})),
$$
$$
Y_\nu^{(2)}(2u;q^2)=\frac{q^{-\nu^2+\frac12}(1-q^2)}
{\pi a_\nu\sqrt{2u}}q^{-\frac12N}
\left[-e^{-i(2\nu-1)\frac\pi4-i\frac{\pi}2n}
e_q^{(2)}(iq^\la)e_q^{(2)}(-iq^{1-\la})+\right.
$$
$$
\left.e^{i(2\nu-1)\frac\pi4+i\frac\pi2n}
e_q^{(2)}(-iq^\la)e_q^{(2)}(iq^{1-\la})\right](1+o(q^{-n})).
$$
\end{predl}
{\it Proof.}  The statement of Proposition follows from (\ref{6.2}) -
(\ref{6.9}). $\Box$

\section{The asymptotic behaviour of the $q^2$-Bessel functions of type 3}
\setcounter{equation}{0}

Consider the coefficient (\ref{5.10}). Obviously
$$
q^{-\nu+\frac12}+q^{\nu-\frac12}\ge2 ~~~\Rightarrow ~~~q^{-\nu+\frac12}-1\ge
1-q^{\nu-\frac12} ~~\Rightarrow
$$
$$
q^{-\frac12}(q^{-\nu+\frac12}-1)\ge q^{\frac12}(1-q^{\nu-\frac12}) ~~\Rightarrow
q^{-\nu}+q^\nu\ge q^{\frac12}+q^{-\frac12} ~~\Rightarrow
$$
$$
(q^\nu+q^{-\nu})q^{l+k+\frac12}\ge(q^{\frac12}+q^{-\frac12})q^{2l+2k+\frac12}
~~\Rightarrow
$$
$$
1-(q^\nu+q^{-\nu})q^{l+k+\frac12}+q^{2l+2k+1}\le1-q^{2l+2k}  ~~\Rightarrow
$$
$$
~~\frac{(1-q^{-\nu+l+k+\frac12})(1-q^{\nu+l+k+\frac12})}{1-q^{2l+2k}}\le1.
$$
Hence it follows from (\ref{1.2}) and (\ref{1.3}) that
$$
\sum_{k=0}^\infty\frac{(q^{-\nu+l+k+\frac12},q)_k(q^{\nu+l+k+\frac12},q)_k}
{(q^{2l+2},q^2)_k(q,q)_k}q^k\le\sum_{k=0}^\infty\frac{q^k}{(q,q)_k}=e_q(q),
$$
$$
\sum_{m=0}^\infty\frac{(q^{-\nu+l+m+\frac12},q)_m(q^{\nu+l+m+\frac12},q)_m}
{(q^{2l+2},q^2)_m(q,q)_m}q^{\frac{m(m+1)}2}\le
\sum_{m=0}^\infty\frac{q^{\frac{m(m+1)}2}}{(q,q)_m}=E_q(q)
$$
and
$$
\sqrt{c_{l-}^{(1)}c_{l-}^{(2)}}\le\sqrt{e_q(q)E_q(q)}
\frac{(q^{-\nu+\frac12},q)_l(q^{\nu+\frac12},q)_l}{(q^2,q^2)_l}q^l.
$$
But this coefficient is not interesting for us because the first summand
of (\ref{5.9}) tends to zero if $|u|\to\infty$.

Consider coefficient (\ref{5.11}). It follows from (\ref{1.1}) and
(\ref{5.3}) that for arbitrary $l\ge0$ and $i=1, 2, ..., k$
$$
1<\frac1{1-q^{l+i}}\le\frac1{1-q},
$$
and for arbitrary $i=1, 2, ..., m$
$$
0<\frac{q^{l+1}}{1-q^{l+i}}<\frac q{1-q}.
$$
So
$$
\sum_{k=0}^\infty\frac{(q^{-\nu+\frac12},q)_k(q^{\nu+\frac12},q)_k}
{(q^2,q^2)_k(q^{l+1},q)_k}q^k=
\sum_{k=0}^\infty\frac{(q^{-\nu+\frac12},q)_k(q^{\nu+\frac12},q)_k}
{(q^2,q^2)_k}(\al q)^k=
$$
\beq{7.1}
\phantom._2\Phi_1(q^{\nu+\frac12},q^{-\nu+\frac12};-q;q,\al q)=
\varphi_1^2(\al)
\eq
for some $\al\in(1,\frac1{1-q})$, and
$$
\sum_{m=0}^\infty\frac{(q^{-\nu+\frac12},q)_m(q^{\nu+\frac12},q)_m}
{(q^2,q^2)_m(q^{l+1},q)_m}q^{\frac12m(m-1)+\be m+m}=
\sum_{m=0}^\infty\frac{(q^{-\nu+\frac12},q)_m(q^{\nu+\frac12},q)_m}
{(q^2,q^2)_m}q^{\frac12m(m-1)}(\be q)^m=
$$
\beq{7.2}
\phantom._2\Phi_2(q^{\nu+\frac12},q^{-\nu+\frac12};-q,0;q,-\be q)=
\varphi_2^2(\be)
\eq
for some $\be\in(0,\frac1{1-q})$.

Hence
\beq{7.3}
\sqrt{c_{l+}^{(1)}c_{l+}^{(2)}}=
\frac{q^{\frac{l(l-1)}4}}{(q,q)_l}\varphi_l(\al)\varphi_2(\be)=
\frac{q^{\frac{l(l-1)}4}}{(q,q)_l}\varphi(\al,\be).
\eq

Now substituting (\ref{7.3}) into (\ref{5.12}) -
(\ref{5.15}) and using (\ref{1.4}) we can write
$$
I_\nu^{(3)}(2u;q^2)=\frac{a_\nu}{\sqrt{2u}}
\left[\sum_{l=1}^\infty(1+ie^{i(\nu+l)\pi})
\sqrt{c_{l-}^{(1)}c_{l-}^{(2)}}u^{-l}+\sum_{l=0}^\infty(1+ie^{i(\nu+l)\pi})
\sqrt{c_{l+}^{(1)}c_{l+}^{(2)}}u^l\right]=
$$
$$
\frac{a_\nu}{\sqrt{2u}}\sum_{l=0}^\infty(1+ie^{i(\nu+l)\pi})
\sqrt{c_{l+}^{(1)}c_{l+}^{(2)}}u^l(1+o(q^{-n}))=
$$
$$
\frac{a_\nu}{\sqrt{2u}}\varphi(\al,\be)(e_q^{(3)}(u)+
ie^{i\nu\pi}e_q^{(3)}(-u))(1+o(q^{-n}))=
$$
$$
\frac{a_\nu}{\sqrt{2u}}q^{-\frac23N-\frac1{24}}
\varphi(\al,\be)\left[\frac{e_q^{(3)}(q^\la)e_q^{(3)}(q^{1-\la})}
{e_q^{(3)}(q^{1-n-\la})}+
ie^{i\nu\pi}\frac{e_q^{(3)}(-q^\la)e_q^{(3)}(-q^{1-\la})}
{e_q^{(3)}(-q^{1-n-\la})}\right](1+o(q^{-n}))
$$
Because $e_q^{(3)}(-q^{1-n-\la})\to1$ if $n\to -\infty$ we have
\beq{7.4}
I_\nu^{(3)}(2u;q^2)=\frac{a_\nu}{\sqrt{2u}}q^{-\frac23N-\frac1{24}}
\varphi(\al,\be)\times
\eq
$$
\times\left[e_q^{(3)}(q^\la)e_q^{(3)}(q^{1-\la})+
ie^{i\nu\pi}e_q^{(3)}(-q^\la)e_q^{(3)}(-q^{1-\la})\right](1+o(q^{-n}))
$$

In the similar way we obtain
\beq{7.5}
K_\nu^{(3)}(2u;q^2)=\frac{q^{-\nu^2+\frac12}(1-q^2)}{2a_\nu\sqrt{2u}}
q^{-\frac23N-\frac1{24}}\varphi(\al,\be)
e_q^{(3)}(-q^\la)e_q^{(3)}(-q^{1-\la})(1+o(q^{-n})),
\eq
\beq{7.6}
J_\nu^{(3)}(2u;q^2)=\frac{2a_\nu}{\sqrt{2z}}q^{-\frac23N-\frac1{24}}
\varphi(\al,\be)e^{-i(2\nu+1)\frac\pi4}\times
\eq
$$
\times\left[-e_q^{(3)}(iq^\la)e_q^{(3)}(-iq^{1-\la})+
ie^{i\nu\pi}e_q^{(3)}(-iq^\la)e_q^{(3)}(iq^{1-\la})\right](1+o(q^{-n})),
$$
\beq{7.7}
Y_\nu^{(3)}(2u;q^2)=-\frac{q^{-\nu^2+\frac12}(1-q^2)}
{a_\nu\pi\sqrt{2u}}q^{-\frac23N-\frac1{24}}
e^{-i(2\nu-1)\frac\pi4}\varphi(\al,\be)\times
\eq
$$
\times\left[e_q^{(3)}(iq^\la)e_q^{(3)}(-iq^{1-\la})+
ie^{i\nu\pi}e_q^{(3)}(-iq^\la)e_q^{(3)}(iq^{1-\la})\right](1+o(q^{-n})),
$$
where $\varphi(\al,\be)=\varphi_1(\al)\varphi_2(\be)$ is determined
be (\ref{7.1}) and (\ref{7.2}).

\bigskip
\small{

}

\end{document}